\documentclass[11pt]{article}
\usepackage{amssymb,latexsym,geometry,cite,enumerate,float}
\newtheorem{theorem}{Theorem}
\newtheorem{proposition}[theorem]{Proposition}
\newtheorem{lemma}[theorem]{Lemma}
\newtheorem{corollary}[theorem]{Corollary}

\newtheorem{conjecture}[theorem]{Conjecture}
\newtheorem{claim}{Claim}

\usepackage{setspace}

\begin{document}

\onehalfspace

\title{Some Bounds on the Zero Forcing Number of a Graph}

\author{Michael Gentner and Dieter Rautenbach}

\date{}

\maketitle

\begin{center}
Institut f\"{u}r Optimierung und Operations Research, 
Universit\"{a}t Ulm, Ulm, Germany,
\{\texttt{michael.gentner, dieter.rautenbach}\}\texttt{@uni-ulm.de}\\[3mm]
\end{center}

\begin{abstract}
A set $Z$ of vertices of a graph $G$ is a zero forcing set of $G$ if 
initially labeling all vertices in $Z$ with $1$ and all remaining vertices of $G$ with $0$,
and then, 
iteratively and as long as possible, 
changing the label of some vertex $u$ from $0$ to $1$
if $u$ is the only neighbor with label $0$ of some vertex with label $1$,
results in the entire vertex set of $G$.
The zero forcing number $Z(G)$, 
defined as the minimum order of a zero forcing set of $G$,
was proposed as an upper bound of the corank of matrices associated with $G$,
and was also considered in connection with quantum physics and logic circuits. 
In view of the computational hardness of the zero forcing number,
upper and lower bounds are of interest.

Refining results of Amos, Caro, Davila, and Pepper, 
we show that $Z(G)\leq \frac{\Delta-2}{\Delta-1}n$ 
for a connected graph $G$ of order $n$ and maximum degree $\Delta$ at least $3$
if and only if $G$ does not belong to $\{ K_{\Delta+1},K_{\Delta,\Delta},K_{\Delta-1,\Delta},G_1,G_2\}$,
where $G_1$ and $G_2$ are two specific graphs of orders $5$ and $7$, respectively.
For a connected graph $G$ of order $n$, maximum degree $3$, and girth at least $5$, 
we show $Z(G)\leq \frac{n}{2}-\Omega\left(\frac{n}{\log n}\right)$.
Using a probabilistic argument,
we show $Z(G)\leq \left(1-\frac{H_r}{r}+o\left(\frac{H_r}{r}\right)\right)n$
for an $r$-regular graph $G$ of order $n$ and girth at least $5$,
where $H_r$ is the $r$-th harmonic number.
Finally, we show 
$Z(G)\geq (g-2)(\delta-2)+2$
for a graph $G$ of girth $g\in \{ 5,6\}$
and minimum degree $\delta$,
which partially confirms a conjecture of Davila and Kenter.
\end{abstract}



\pagebreak

\section{Introduction}

We consider graphs that are finite, simple, and undirected, and use standard terminology.

Let $G$ be a graph.
For a set $Z$ of vertices of $G$, let ${\cal F}(Z)$ be the maximal set of vertices of $G$ that arises from $Z$
by iteratively adding vertices that are the unique neighbor outside the current set of some vertex inside the current set.
Equivalently, 
\begin{itemize}
\item $|N_G(w)\setminus {\cal F}(Z)|\not=1$ for every vertex $w$ in ${\cal F}(Z)$, and,
\item the elements of ${\cal F}(Z)\setminus Z$ have a linear order $u_1,\ldots,u_k$ such that
for every index $i$ in $\{ 1,\ldots,k\}$, 
there is some vertex $v_i$ in $Z\cup \{ u_j:1\leq j\leq i-1\}$ such that $u_i$ is the only neighbor of $v_i$ in $\{ u_j:i\leq j\leq k\}$.
\end{itemize}
In the latter case, we say that $v_i$ {\it forces} $u_i$ for $i\in \{ 1,\ldots,k\}$, and denote this by $v_i\to u_i$.
The sequence $v_1\to u_1, v_2\to u_2,\ldots,v_k\to u_k$ is called a {\it forcing sequence} for $Z$.

The set $Z$ is a {\it zero forcing set} of $G$ if ${\cal F}(Z)$ equals the vertex set $V(G)$ of $G$.
The {\it zero forcing number $Z(G)$} of $G$ is the minimum order of a zero forcing set of $G$.
The zero forcing number was proposed by the AIM Minimum Rank - Special Graphs Work Group \cite{aim,hv} 
as an upper bound on the corank of matrices associated with a given graph.
Independently, it was considered in connection with 
quantum physics \cite{bg,bm,s}
as well as logic circuits \cite{bghsy}.
It has already been studied in a number of papers \cite{acdp,cdky,ehhlr,hhkmwy,m,r,r2,dk,td}
and is computationally hard \cite{a,fmy}.

\medskip

\noindent In the present paper we establish some upper and lower bounds on the zero forcing number.

For a connected graph $G$ of order $n$ and maximum degree $\Delta$ at least $2$,
Amos et al. \cite{acdp} prove
\begin{eqnarray}
Z(G)&\leq &\frac{\Delta}{\Delta+1}n\mbox{ and}\label{eb1}\\
Z(G)&\leq &\frac{\Delta-2}{\Delta-1}n+\frac{2}{\Delta+1}.\label{eb2}
\end{eqnarray}
It was shown that the only extremal graph for (\ref{eb1}) is the complete graph $K_{\Delta+1}$ of order $\Delta+1$ \cite{gprs},
and that the only extremal graphs for (\ref{eb2}) are $K_{\Delta+1}$,
the complete bipartite graph $K_{\Delta,\Delta}$ with partite sets of order $\Delta$, and
the cycle $C_{n}$ \cite{gprs,lwz}.

We characterize the graphs for which the additive term $\frac{2}{\Delta+1}$ in (\ref{eb2}) is not needed.
In fact, we believe that (\ref{eb2}) can be improved considerably,
and, in particular, pose the following conjecture.

\begin{conjecture}\label{conjecture2}
If $G$ is a connected graph of order $n$ and maximum degree $3$, then $Z(G)\leq \frac{1}{3}n+2$.
\end{conjecture}
As a contribution towards this conjecture,
we prove $Z(G)\leq \frac{n}{2}-\Omega\left(\frac{n}{\log n}\right)$
for a connected graph $G$ of order $n$, maximum degree $3$, and girth at least $5$.
We present a probabilistic upper bound on the zero forcing number and discuss some of its consequences.

In \cite{dk} Davila and Kenter conjecture that the lower bound
\begin{eqnarray}
Z(G)&\geq &(g-2)(\delta-2)+2\label{eb3}
\end{eqnarray}
for every graph $G$ of girth $g$ at least $3$ and minimum degree $\delta$ at least $2$.
They observe that for $g>6$ and sufficiently large $\delta$,
the conjecture follows by combining results from \cite{aim2} and \cite{cs}.
For $g=4$, that is, for triangle-free graphs, it was shown in \cite{gprs}.
Here, we prove the conjecture for $g\in \{ 5,6\}$.

\section{Results}

We begin with a simple consequence of (\ref{eb2}).

\begin{proposition}\label{proposition1}
If $G$ is a connected graph of order $n$ and maximum degree $\Delta$ at least $3$ that is distinct from $K_{\Delta+1}$, then 
$$Z(G)\leq \frac{\Delta-1}{\Delta}n.$$
\end{proposition}
{\it Proof:} If $n\geq 2\Delta$, then $Z(G)\stackrel{(\ref{eb2})}{\leq} \frac{(\Delta-2)n+2}{\Delta-1}\leq \frac{(\Delta-1)n}{\Delta}$.
Now, let $n<2\Delta$.
Since $G$ is not complete, it contains an induced path $uvw$ of order $3$.
Since the set $V(G)\setminus \{ v,w\}$ is a zero-forcing set of $G$,
we obtain $Z(G)\leq n-2\leq \frac{(\Delta-1)n}{\Delta}$,
which completes the proof. $\Box$

\medskip

\noindent Our next goal is to characterize the graphs for which the additive term in (\ref{eb2}) is not needed.

The following lemma is implicit in the greedy argument in \cite{cp}.

\begin{lemma}\label{lemma1}
Let $G$ be a connected graph of order $n$ and maximum degree $\Delta$ at least $3$.

If there is some set $Z_0$ of vertices of $G$ such that 
$|{\cal F}(Z_0)|\geq \frac{\Delta-1}{\Delta-2}|Z_0|$,
and ${\cal F}(Z_0)$ induces a subgraph of $G$ without isolated vertices,
then 
$Z(G)\leq \frac{\Delta-2}{\Delta-1}n$.
\end{lemma}
{\it Proof:} If ${\cal F}(Z_0)=V(G)$, then $Z_0$ is a zero forcing set of $G$, and, hence, 
$Z(G)\leq |Z_0|\leq \frac{\Delta-2}{\Delta-1}|{\cal F}(Z_0)|=\frac{\Delta-2}{\Delta-1}n$.
Therefore, we may assume that $Z_i$ is a set of vertices of $G$ for some non-negative integer $i$
such that $|{\cal F}(Z_i)|\geq \frac{\Delta-1}{\Delta-2}|Z_i|$,
the set ${\cal F}(Z_i)$ induces a subgraph of $G$ without isolated vertices,
and ${\cal F}(Z_i)$ is a proper subset of $V(G)$.
Because $G$ is connected, there is a vertex $v$ in ${\cal F}(Z_i)$ 
that has at least one neighbor in $V(G)\setminus {\cal F}(Z_i)$ as well as at least one neighbor in ${\cal F}(Z_i)$.
Let $Z_{i+1}$ arise from $Z_i$ by adding to $Z_i$ all but exactly one neighbor of $v$ in $V(G)\setminus {\cal F}(Z_i)$.
Note that $|Z_{i+1}|=|Z_i|+|N_G(v)\setminus {\cal F}(Z_i)|-1$.
Since $N_G(v)\subseteq {\cal F}(Z_{i+1})$, we obtain $|{\cal F}(Z_{i+1})|\geq {\cal F}(Z_i)+|N_G(v)\setminus {\cal F}(Z_i)|$.
Since $|N_G(v)\setminus {\cal F}(Z_i)|\leq \Delta-1$, this implies 
$|{\cal F}(Z_{i+1})|\geq \frac{\Delta-1}{\Delta-2}|Z_{i+1}|$.
Furthermore, by construction, the set ${\cal F}(Z_{i+1})$ induces a subgraph of $G$ without isolated vertices.
Repeating this extension as long as ${\cal F}(Z_i)$ is a proper subset of $V(G)$, we obtain a zero forcing set $Z$
of $G$ with $|Z|\leq \frac{\Delta-2}{\Delta-1}|{\cal F}(Z)|=\frac{\Delta-2}{\Delta-1}n$,
which completes the proof. $\Box$

\begin{figure}[H]
\begin{center}
$\mbox{}$\hfill
\unitlength 1mm 
\linethickness{0.4pt}
\ifx\plotpoint\undefined\newsavebox{\plotpoint}\fi 
\begin{picture}(21,18)(0,0)
\put(0,0){\circle*{2}}
\put(15,0){\circle*{2}}
\put(10,10){\circle*{2}}
\put(20,10){\circle*{2}}
\put(0,10){\circle*{2}}
\put(0,10){\line(1,0){20}}
\put(20,10){\line(-1,-2){5}}
\put(15,0){\line(-1,2){5}}
\put(0,10){\line(0,-1){10}}
\put(0,0){\line(1,0){15}}
\qbezier(20,10)(10,18)(0,10)
\end{picture}\hfill
\unitlength 1mm 
\linethickness{0.4pt}
\ifx\plotpoint\undefined\newsavebox{\plotpoint}\fi 
\begin{picture}(26,28)(0,0)
\put(5,10){\circle*{2}}
\put(20,10){\circle*{2}}
\put(15,20){\circle*{2}}
\put(25,20){\circle*{2}}
\put(5,20){\circle*{2}}
\put(5,20){\line(1,0){20}}
\put(25,20){\line(-1,-2){5}}
\put(20,10){\line(-1,2){5}}
\put(5,20){\line(0,-1){10}}
\put(5,10){\line(1,0){15}}
\qbezier(25,20)(15,28)(5,20)
\put(20,0){\circle*{2}}
\put(5,0){\circle*{2}}
\put(15,20){\line(1,-4){5}}
\put(20,0){\line(1,4){5}}
\put(20,10){\line(-3,-2){15}}
\put(5,0){\line(1,0){15}}
\put(20,0){\line(-3,2){15}}
\put(5,10){\line(0,-1){10}}
\qbezier(5,20)(-3,10)(5,0)
\end{picture}\hfill$\mbox{}$
\end{center}
\caption{The two specific graphs $G_1$ and $G_2$.}\label{fig0}
\end{figure}
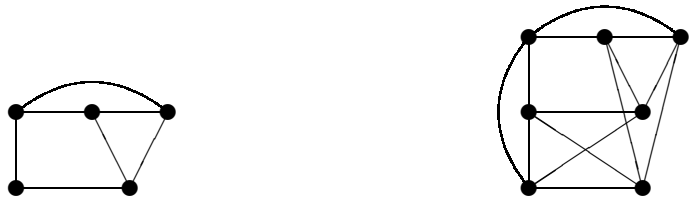

\begin{theorem}\label{theorem4}
If $G$ is a connected graph of order $n$ and maximum degree $\Delta$ at least $3$, then 
\begin{eqnarray}\label{eb3}
Z(G)\leq \frac{\Delta-2}{\Delta-1}n
\end{eqnarray}
if and only if $G\not\in \{ K_{\Delta+1},K_{\Delta,\Delta},K_{\Delta-1,\Delta},G_1,G_2\}$,
where $G_1$ and $G_2$ are the two specific graphs illustrated in Figure \ref{fig0}.
\end{theorem}
{\it Proof:} The necessity follows easily using 
$Z(K_{\Delta+1})=\Delta$,
$Z(K_{\Delta,\Delta})=2\Delta-2$,
$Z(K_{\Delta-1,\Delta})=2\Delta-3$,
$Z(G_1)=3$, and
$Z(G_2)=5$.
We proceed to the proof of the sufficiency.
Therefore, let $G\not\in \{ K_{\Delta+1},K_{\Delta,\Delta},K_{\Delta-1,\Delta},G_1,G_2\}$ be as in the statement.
In order to derive (\ref{eb3}) using Lemma \ref{lemma1},
it suffices to exhibit a set $Z_0$ of vertices of $G$ such that
\begin{eqnarray}\label{el1}
\mbox{$\frac{|{\cal F}(Z_0)|}{|Z_0|}\geq \frac{\Delta-1}{\Delta-2}$, and ${\cal F}(Z_0)$ induces a subgraph of $G$ without isolated vertices.}
\end{eqnarray}
Therefore, suppose that such a set does not exist.

If $G$ has a vertex $v$ of degree $d_G(v)$ at most $\Delta-2$, and $u$ is a neighbor of $v$,
then let $Z_0=N_G[v]\setminus \{ u\}$. 
Since $|Z_0|=d_G(v)$ and $u\in {\cal F}(Z_0)$, 
we obtain $\frac{|{\cal F}(Z_0)|}{|Z_0|}\geq \frac{d_G(v)+1}{d_G(v)}\geq \frac{\Delta-1}{\Delta-2}$,
that is, the set $Z_0$ satisfies (\ref{el1}), which is a contradiction.
Hence, we may assume that $G$ has minimum degree at least $\Delta-1$.

Since $\Delta-1\geq 2$, the graph $G$ is not a tree.
Let $C:v_1\ldots v_gv_1$ be a shortest cycle in $G$.
We consider three cases depending on the girth $g$ of $G$.

\medskip

\noindent {\bf Case 1} $g\geq 5$.

\medskip

\noindent Since $G$ has girth at least $5$, no vertex in $V(G)\setminus V(C)$ has more than one neighbor on $C$.
If all vertices on $C$ have degree at least $3$, 
then let $u_i$ be a neighbor of $v_i$ in $V(G)\setminus V(C)$ for every $i\in \{ 1,\ldots,g\}$.
Let $Z_0=\bigcup\limits_{i=1}^g N_G[v_i]\setminus \{ u_i\}$.
Since $|Z_0|\leq (\Delta-2)g$ and $u_1,\ldots,u_g\in {\cal F}(Z_0)$, 
we obtain $\frac{|{\cal F}(Z_0)|}{|Z_0|}\geq \frac{|Z_0|+g}{|Z_0|}\geq \frac{\Delta-1}{\Delta-2}$,
that is, the set $Z_0$ satisfies (\ref{el1}), which is a contradiction.
Hence, we may assume that $C$ contains a vertex of degree $2$.
Since $G$ has minimum degree at least $\Delta-1\geq 2$, this implies $\Delta=3$.

Let $1\leq i_1<i_2<\ldots<i_p\leq g$ be such that $\{ v_{i_j}:1\leq j\leq p\}$ is the set of vertices of degree $3$ on $C$.
Since $G$ is connected and has maximum degree $3$, we obtain that $p$ is at least $1$.
Possibly renaming vertices, we may assume that $i_p=g$.
Similarly as above, for $j\in \{ 1,\ldots,p\}$, let $u_{i_j}$ be the neighbor of $v_{i_j}$ in $V(G)\setminus V(C)$.

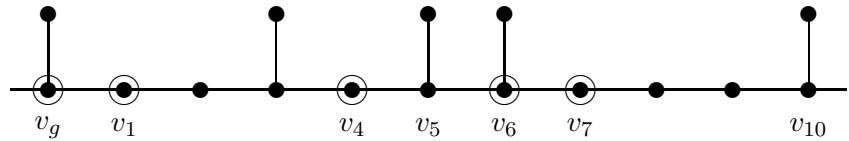
\begin{figure}[H]
\begin{center}
\unitlength 1mm 
\linethickness{0.4pt}
\ifx\plotpoint\undefined\newsavebox{\plotpoint}\fi 
\begin{picture}(110,16)(0,0)
\put(5,5){\circle*{2}}
\put(5,15){\circle*{2}}
\put(15,5){\circle*{2}}
\put(25,5){\circle*{2}}
\put(35,5){\circle*{2}}
\put(35,15){\circle*{2}}
\put(45,5){\circle*{2}}
\put(55,5){\circle*{2}}
\put(55,15){\circle*{2}}
\put(65,5){\circle*{2}}
\put(65,15){\circle*{2}}
\put(75,5){\circle*{2}}
\put(85,5){\circle*{2}}
\put(95,5){\circle*{2}}
\put(105,5){\circle*{2}}
\put(105,15){\circle*{2}}
\put(0,5){\line(1,0){110}}
\put(5,15){\line(0,-1){10}}
\put(35,15){\line(0,-1){10}}
\put(55,15){\line(0,-1){10}}
\put(65,15){\line(0,-1){10}}
\put(105,15){\line(0,-1){10}}
\put(5,5){\circle{4}}
\put(15,5){\circle{4}}
\put(45,5){\circle{4}}
\put(65,5){\circle{4}}
\put(75,5){\circle{4}}
\put(5,0){\makebox(0,0)[cc]{$v_g$}}
\put(15,0){\makebox(0,0)[cc]{$v_1$}}
\put(45,0){\makebox(0,0)[cc]{$v_4$}}
\put(55,0){\makebox(0,0)[cc]{$v_5$}}
\put(65,0){\makebox(0,0)[cc]{$v_6$}}
\put(75,0){\makebox(0,0)[cc]{$v_7$}}
\put(105,0){\makebox(0,0)[cc]{$v_{10}$}}
\end{picture}
\caption{A section of $C$. The vertices in $\{ v_1,v_4,v_6,v_7\}$ belong to $Z_0$ because the vertices in $\{ v_g,v_3,v_5,v_6\}$ have degree $3$. The vertex $v_g$ belongs to $Z_0$ regardless of the degree of $v_{g-1}$.
Note that 
$v_1\to v_2$, 
$v_2\to v_3$, 
$v_3\to u_3$, 
$v_4\to v_5$, 
$v_5\to u_5$, 
$v_6\to u_6$, 
$v_7\to v_8$, 
$v_8\to v_9$, 
$v_9\to v_{10}$, and
$v_{10}\to u_{10}$.}\label{fig1}
\end{center}
\end{figure}

If $p\leq g-2$, then let $Z_0=\{ v_g\}\cup \{ v_{i_j+1}:1\leq j\leq p\}$, where $v_{g+1}=v_1$.
See Figure \ref{fig1} for an illustration.

\noindent Since $|Z_0|\leq p+1$ and $V(C)\cup \{ u_{i_j}:1\leq j\leq p\}\subseteq {\cal F}(Z_0)$,
we obtain 
$\frac{|{\cal F}(Z_0)|}{|Z_0|}\geq \frac{g+p}{p+1}\geq 2=\frac{\Delta-1}{\Delta-2}$,
that is, the set $Z_0$ satisfies (\ref{el1}), which is a contradiction.
Hence, we may assume that $p=g-1$, that is, $C$ contains exactly one vertex, say $v_1$, of degree $2$.
Let $Z_0=V(C)\setminus \{ v_2\}$.
Since $|Z_0|\leq g-1$ and $V(C)\cup \{ u_{i_j}:1\leq j\leq p\}\subseteq {\cal F}(Z_0)$,
we obtain 
$\frac{|{\cal F}(Z_0)|}{|Z_0|}\geq \frac{2g-1}{g-1}\geq 2=\frac{\Delta-1}{\Delta-2}$,
that is, the set $Z_0$ satisfies (\ref{el1}), which is a contradiction.
This completes the proof in this case.

\medskip

\noindent {\bf Case 2} $g=4$.

\medskip

\noindent First, we assume that $d_G(v_1)=2$.
As noted above, this implies $\Delta=3$.
If $d_G(v_2)=2$, then $Z_0=\{ v_1,v_2\}$ satisfies (\ref{el1}), which is a contradiction.
Hence, by symmetry, we may assume that $d_G(v_2)=d_G(v_4)=3$.
Let $Z_0=\{ v_1,v_2,v_3\}$.
If $v_1$ and $v_3$ are the only common neighbors of $v_2$ and $v_4$, 
then ${\cal F}(Z_0)$ contains $v_4$ as well as the two neighbors of $v_2$ and $v_4$ that do not lie on $C$.
Hence, $|{\cal F}(Z_0)|\geq 6$,
that is, the set $Z_0$ satisfies (\ref{el1}), which is a contradiction.
Hence, we may assume that $N_G(v_2)=N_G(v_4)$.
Since $G\not=K_{2,3}$,
we may assume, by symmetry, that $d_G(v_3)=3$.
Since ${\cal F}(Z_0)$ contains $N_G[v_2]$ and the neighbor of $v_3$ that does not lie on $C$,
we obtain $|{\cal F}(Z_0)|\geq 6$,
that is, the set $Z_0$ satisfies (\ref{el1}), which is a contradiction.
Hence, we may assume, by symmetry, that $G$ contains no cycle of length $4$ that contains a vertex of degree $2$.
Since $G$ is a shortest cycle, it is induced.
For $i\in \{ 1,2,3,4\}$, let $u_i$ be a neighbor of $v_i$ that does not lie on $C$.

\medskip

\noindent Next, we assume that $N_G(v_1)\not\subseteq N_G(v_3)$ and $N_G(v_2)\not\subseteq N_G(v_4)$.
We may assume that $u_1\in N_G(v_1)\setminus N_G(v_3)$ and $u_2\in N_G(v_2)\setminus N_G(v_4)$,
which implies that $u_1$, $u_2$, $u_3$, and $u_4$ are four distinct vertices.
Let $Z_0=(N_G[v_1]\cup N_G[v_2]\cup N_G[v_3]\cup N_G[v_4])\setminus \{ u_1,u_2,u_3,u_4\}$.
Clearly, $|Z_0|\leq 4(\Delta-2)$.
Since 
$v_1\to u_1$, 
$v_2\to u_2$, 
$v_3\to u_3$, and 
$v_4\to u_4$, 
we obtain $u_1,u_2,u_3,u_4\in {\cal F}(Z_0)$, and, hence,
$\frac{|{\cal F}(Z_0)|}{|Z_0|}\geq \frac{|Z_0|+4}{|Z_0|}\geq \frac{4(\Delta-2)+4}{4(\Delta-2)}=\frac{\Delta-1}{\Delta-2}$,
that is, the set $Z_0$ satisfies (\ref{el1}), which is a contradiction.
Hence, we may assume, by symmetry, that $N_G(u_2)=N_G(u_4)$.

\medskip

\noindent Next, we assume that $N_G(v_1)\not\subseteq N_G(v_3)$.
Again, let $u_1\in N_G(v_1)\setminus N_G(v_3)$.
If $|N_G(v_1)\cup N_G(v_3)|\leq 2\Delta-3$, 
then let 
$Z_0=(N_G[v_1]\cup N_G[v_2]\cup N_G[v_3])\setminus \{ u_1,u_2,u_3\}$.
We obtain $|Z_0|\leq |N_G(v_1)\cup N_G(v_3)|+|N_G(v_2)|-|\{ u_1,u_2,u_3\}|
\leq 2\Delta-3+\Delta-3=3(\Delta-2)$.
Since 
$v_1\to u_1$, 
$v_2\to u_2$, and
$v_3\to u_3$, 
we obtain $u_1,u_2,u_3\in {\cal F}(Z_0)$, and, hence,
$\frac{|{\cal F}(Z_0)|}{|Z_0|}\geq \frac{|Z_0|+3}{|Z_0|}\geq \frac{3(\Delta-2)+3}{3(\Delta-2)}=\frac{\Delta-1}{\Delta-2}$,
that is, the set $Z_0$ satisfies (\ref{el1}), which is a contradiction.
Hence, we may assume $|N_G(v_1)\cup N_G(v_3)|\geq 2\Delta-2$,
which implies that $v_1$ and $v_3$ both have degree $\Delta$, and do not have a common neighbor apart from $v_2$ and $v_4$.
By symmetry, this implies that every vertex in $N_G(u_2)$ has degree $\Delta$,
and that every two vertices in $N_G(u_2)$ do not have a common neighbor apart from $v_2$ and $v_4$.
Let $w_2\in N_G(u_2)\setminus \{ v_2,v_4\}$,
and let 
$Z_0=(N_G[v_1]\cup N_G[v_2]\cup N_G[v_3]\cup N_G[u_2])\setminus \{ u_1,u_2,u_3,w_2\}$.
We obtain $|Z_0|\leq 4(\Delta-2)$.
Since 
$v_2\to u_2$,
$v_1\to u_1$,
$v_3\to u_3$, and
$u_2\to w_2$, 
we obtain $u_1,u_2,u_3,w_2\in {\cal F}(Z_0)$, and, hence,
$\frac{|{\cal F}(Z_0)|}{|Z_0|}\geq \frac{\Delta-1}{\Delta-2}$,
that is, the set $Z_0$ satisfies (\ref{el1}), which is a contradiction.
Hence, we may assume that $N_G(v_1)=N_G(v_3)$.

\medskip

\noindent If some vertex $v_4'\in N_G(v_1)$ is not adjacent to some vertex in $N_G(v_2)$, 
then one of the previous cases applies to the cycle $v_1v_2v_3v_4'v_1$.
Hence, we may assume that all vertices in $N_G(v_1)$ are adjacent to all vertices in $N_G(v_2)$,
which implies that $G$ contains a complete bipartite subgraph $H$ with partite sets $N_G(v_1)$ and $N_G(v_2)$.
If $N_G(v)\not=N_G(w)$ for two vertices $v$ and $w$ that both either belong to $N_G(v_1)$ or to $N_G(v_2)$,
then some previous case applies to a cycle of length $4$ containing these two vertices.
This implies that $G$ equals $H$, and, hence, $Z(G)=n-2$.
Since $G\not\in \{ K_{\Delta,\Delta},K_{\Delta-1,\Delta}\}$,
we obtain $n\leq 2\Delta-2$ and (\ref{eb3}) follows,
which completes the proof in this case.

\medskip

\noindent {\bf Case 3} $g=3$.

\medskip

\noindent First, we assume that $d_G(v_1)=2$.
Again, this implies $\Delta=3$.
Since $G$ is connected and has maximum degree $3$, we may assume that $d_G(v_2)=3$.
This implies that the set $Z_0=\{ v_1,v_2\}$ satisfies (\ref{el1}), which is a contradiction.
Hence, we may assume, by symmetry, that $G$ contains no triangle that contains a vertex of degree $2$.
For $i\in \{ 1,2,3\}$, let $u_i$ be a neighbor of $v_i$ that does not lie on $C$.

\medskip

\noindent Next, we assume that $N_G(v_1)\not\subseteq N_G(v_2)\cup N_G(v_3)$ and $N_G(v_2)\not\subseteq N_G(v_3)$.
We may assume that $u_1\in N_G(v_1)\setminus (N_G(v_2)\cup N_G(v_3))$ and $u_2\in N_G(v_2)\setminus N_G(v_3)$.
For $Z_0=(N_G[v_1]\cup N_G[v_2]\cup N_G[v_3])\setminus \{ u_1,u_2,u_3\}$,
we obtain $|Z_0|\leq 3(\Delta-2)$.
Since 
$v_3\to u_3$,
$v_2\to u_2$, and
$v_1\to u_1$, 
we obtain $u_1,u_2,u_3\in {\cal F}(Z_0)$, and, hence,
$\frac{|{\cal F}(Z_0)|}{|Z_0|}\geq \frac{\Delta-1}{\Delta-2}$,
that is, the set $Z_0$ satisfies (\ref{el1}), which is a contradiction.

\medskip

\noindent Next, we assume that $N_G(v_1)\not\subseteq N_G(v_2)\cup N_G(v_3)$ and $N_G(v_2)=N_G(v_3)$.
If $|N_G(v_1)\cup N_G(v_2)|\leq 2\Delta-2$, 
then let $Z_0=(N_G[v_1]\cup N_G[v_2])\setminus \{ u_1,u_2\}$.
Note that $|Z_0|\leq 2(\Delta-2)$.
Since
$v_2\to u_2$
and
$v_1\to u_1$,
we obtain $u_1,u_2\in {\cal F}(Z_0)$, and, hence,
$\frac{|{\cal F}(Z_0)|}{|Z_0|}\geq \frac{\Delta-1}{\Delta-2}$,
that is, the set $Z_0$ satisfies (\ref{el1}), which is a contradiction.
Hence, $|N_G(v_1)\cup N_G(v_2)|\geq 2\Delta-1$, 
which implies that $v_1$, $v_2$, and $v_3$ all have degree $\Delta$,
and that $v_3$ is the only common neighbor of $v_1$ and $v_2$.
Let $A=N_G(v_1)\setminus \{ v_2,v_3\}$ and $B=N_G(v_2)\setminus \{ v_1,v_3\}$.
Note that $|A|=|B|=\Delta-2$.
If some vertex $u_1'$ in $A$ is not adjacent to some vertex $u_2'$ in $B$,
then let $Z_0=(N_G(v_1)\cup N_G(v_2)\cup N_G(u_2'))\setminus \{ u_1',v_1,v_2\}$.
Note that $|Z_0|\leq 3(\Delta-2)$.
Since
$u_2'\to v_2$,
$v_2\to v_1$, and 
$v_1\to u_1'$,
we obtain $v_1,v_2,u_1'\in {\cal F}(Z_0)$, and, hence,
$\frac{|{\cal F}(Z_0)|}{|Z_0|}\geq \frac{\Delta-1}{\Delta-2}$,
that is, the set $Z_0$ satisfies (\ref{el1}), which is a contradiction.
Hence, every vertex in $A$ is adjacent to every vertex in $B$.
Note that $N_G(u)=\{ v_2,v_3\}\cup A$ for every vertex $u$ in $B$,
and that every vertex in $A$ has at most one neighbor outside of $\{ v_1\}\cup B$.

If some vertex $u_1'$ in $A$ has a neighbor $w_1$ outside of $\{ v_1,v_2,v_3\}\cup A\cup B$,
then let $Z_0=(N_G(v_1)\cup N_G(v_2))\setminus \{ u_1',u_2\}$.
Note that $|Z_0|=2\Delta-3$.
Since
$v_2\to u_2$,
$v_1\to u_1'$, and 
$u_1'\to w_1$,
we obtain $u_1',u_2,w_1\in {\cal F}(Z_0)$, and, hence,
$\frac{|{\cal F}(Z_0)|}{|Z_0|}\geq \frac{|Z_0|+3}{|Z_0|}=\frac{2\Delta}{2\Delta-3}\geq \frac{\Delta-1}{\Delta-2}$,
that is, the set $Z_0$ satisfies (\ref{el1}), which is a contradiction.
Hence, no vertex in $A$ has a neighbor outside of $\{ v_1,v_2,v_3\}\cup A\cup B$.
Note that $A$ induces a subgraph of $G$ of maximum degree at most $1$.

If $A$ contains two vertices $u_1'$ and $u_1''$ that are not adjacent,
then let $Z_0=(N_G(v_1)\cup N_G(v_2))\setminus \{ u_1'',v_2,u_2\}$.
Note that $|Z_0|=2\Delta-4$.
Since
$u_1'\to u_2$,
$v_3\to v_2$, and 
$v_1\to u_1''$,
we obtain $u_1'',v_2,u_2\in {\cal F}(Z_0)$, and, hence,
$\frac{|{\cal F}(Z_0)|}{|Z_0|}\geq \frac{|Z_0|+3}{|Z_0|}=\frac{2\Delta-1}{2\Delta-4}\geq \frac{\Delta-1}{\Delta-2}$,
that is, the set $Z_0$ satisfies (\ref{el1}), which is a contradiction.
Hence, every two vertices in $A$ are adjacent.

Since $G$ has maximum degree $\Delta$, and every vertex in $A$ has degree
$|\{ v_1\}|+(|A|-1)+|B|=2\Delta-4$, we obtain $\Delta\leq 4$,
which implies the contradiction that $G$ is either $G_1$ or $G_2$.
Hence, we may assume, by symmetry, that 
$N_G(v_i)\subseteq N_G(v_j)\cup N_G(v_k)$ for $\{ i,j,k\}=\{ 1,2,3\}$.
Note that this implies 
$|N_G[v_1]\cup N_G[v_2]\cup N_G[v_3]|\leq \frac{3}{2}(\Delta-2)+3$.

\medskip

\noindent Since $G$ is not $K_{\Delta+1}$, we may assume that $N_G(v_1)\not\subseteq N_G(v_2)$, and that $\Delta\geq 4$.
We may assume that $u_1\in N_G(v_1)\setminus N_G(v_2)$.
Let $Z_0=(N_G(v_1)\cup N_G(v_2))\setminus \{ u_1,u_2\}$.
Note that $|Z_0|\leq \frac{3}{2}(\Delta-2)+1$.
Since 
$v_2\to u_2$ and
$v_1\to u_1$,
we obtain $u_1,u_2\in {\cal F}(Z_0)$, and, hence,
$\frac{|{\cal F}(Z_0)|}{|Z_0|}\geq \frac{|Z_0|+2}{|Z_0|}
\geq 
\frac{\frac{3}{2}(\Delta-2)+3}{\frac{3}{2}(\Delta-2)+1}
\geq \frac{\Delta-1}{\Delta-2}$,
that is, the set $Z_0$ satisfies (\ref{el1}), which is a contradiction.
This completes the proof.
$\Box$

\medskip

\noindent While our Conjecture \ref{conjecture2} remains widely open, 
we are able to improve (\ref{eb2}) at least by some lower order term
for subcubic graphs of girth at least $5$.

\begin{theorem}\label{theorem3}
If $G$ is a connected graph of order $n$, maximum degree $3$, and girth at least $5$, then 
$$Z(G)\leq \frac{n}{2}-\frac{n}{24\log_2(n)+6}+2.$$
\end{theorem}
{\it Proof:} Let $G$ be as in the statement.
We begin with an extension statement similar to Lemma \ref{lemma1}.

\begin{claim}\label{claim1}
Let $Z$ be a set of vertices of $G$.
Let $F={\cal F}(Z)$ and $R=V(G)\setminus F$.

If $F$ induces a connected subgraph of $G$ of order at least $3$, and $R$ contains a vertex of degree at least $2$, 
then there is a set $Z'$ of vertices of $G$ with 
\begin{enumerate}[(i)]
\item $|Z'\setminus Z|\leq 2\log_2(n)$,  
\item $|{\cal F}(Z')\setminus F|\geq 2|Z'\setminus Z|+1$,
\item $Z\subseteq Z'$, and ${\cal F}(Z')$ induces a connected subgraph of $G$.
\end{enumerate}
\end{claim}
{\it Proof of Claim \ref{claim1}:}
Note that a vertex in $F$ with a neighbor in $R$ has exactly one neighbor in $F$ and two neighbors in $R$,
in particular, such a vertex has degree $3$.

A subgraph $H$ of $G$ is an {\it extension subgraph} if it is of one of the following types:
\begin{enumerate}[Type a:]
\item A path $P:v_0\ldots v_k$ with $v_0\in F$, $v_1,\ldots,v_k\in R$, and $d_G(v_k)=2$.
\item A path $P:v_0\ldots v_k$ with $v_0\in F$, $v_1,\ldots,v_k\in R$, $d_G(v_k)=1$, and $k\geq 2$.
\item A path $P:v_0\ldots v_k$ with $v_0,v_k\in F$, $v_1,\ldots,v_{k-1}\in R$, and $k\geq 2$.
\item A cycle $C:u_1\ldots u_\ell u_1$ with $u_1\in F$, $u_2,\ldots,u_\ell \in R$.
\item The union of a path $P:v_0\ldots v_k$ and a cycle $C:u_1\ldots u_\ell u_1$
with $v_0\in F$, $v_1,\ldots,v_k,u_1,\ldots,u_\ell \in R$, $v_k=u_1$, and $V(P)\cap V(C)=\{ u_1\}$.
\end{enumerate}
Whenever we refer to some extension subgraph,
we use the notation introduced above.

First, we show the existence of a small extension subgraph.
Therefore, suppose that $G$ does not contain an extension subgraph of order at most $2\log_2(n)+1$.
Since $G$ is connected, and $R$ contains a vertex of degree at least $2$, 
there is a vertex $v$ in $F$ that has a neighbor $u$ in $R$ such that $u$ has degree at least $2$.
Since there is no extension subgraph of order at most $2\log_2(n)+1$,
the vertex $u$ is the root of a perfect binary subtree $T$ of $G$ of height $\lfloor\log_2(n)\rfloor$ with $V(T)\subseteq R$.
Since $v$ has a neighbor in $F$, we obtain the contradiction $n\geq 2+n(T)=2+2^{\lfloor\log_2(n)\rfloor+1}-1>n$. 

\medskip

\noindent Let $H$ be an extension subgraph
such that the order $n(H)$ of $H$ is as small as possible,
and, subject to this first condition,
the number of vertices of $H$ in $R$ is as small as possible.

As shown above, $n(H)\leq 2\log_2(n)+1$.

Since $G$ has girth at least $5$, and the set $F$ contains more than two vertices, 
the choice of $H$ easily implies that 
\begin{itemize}
\item $H$ is an induced subgraph of $G$,
\item no vertex in $R\setminus V(H)$ is adjacent to two vertices of $H$, 
\item $V(H)\cap R$ contains a vertex of degree less than $3$ only if $H$ has Type a or Type b, 
in which case $v_k$ is the only such vertex,
and
\item every vertex $v$ in $V(H)\cap R$ with $d_H(u)=2$ and $d_G(u)=3$ has a neighbor $p(v)$ in $R\setminus V(H)$.
\end{itemize}
The violation of any of these conditions leads to an extension subgraph of smaller order
or of the same order but less vertices in $R$.
As observed above, every vertex $v$ in $V(H)\cap F$ has exactly two neighbors in $R$,
and if only one of these two neighbors belongs to $H$, then we denote the other neighbor by $p(v)$.

\medskip

\noindent Now, we consider the different types.

First, assume that $H$ has Type a).
Let $u$ be the neighbor of $v_k$ distinct from $v_{k-1}$.
If $u\in F$, then the choice of $H$ implies $k=1$.
Let $Z'=Z\cup \{ v_k\}$, and let $p(u)$ be the neighbor of $u$ in $R$ distinct from $v_k$.
Since $|Z'\setminus Z|=1$, we obtain (i).
Since $p(u),p(v_0)\in {\cal F}(Z')\setminus F$ and $Z'\setminus Z\subseteq {\cal F}(Z')\setminus F$, 
we obtain $|{\cal F}(Z')\setminus F|\geq 2|Z'\setminus Z|+1$, and, hence, (ii).
If $u\not\in F$, then let $Z'=Z\cup \{ p(v_i):0\leq i\leq k-1\}$.
Since $|Z'\setminus Z|=k=n(H)-1\leq 2\log_2(n)$, we obtain (i).
Since $v_1,\ldots,v_k,u\in {\cal F}(Z')\setminus F$ and $Z'\setminus Z\subseteq {\cal F}(Z')\setminus F$, 
we obtain $|{\cal F}(Z')\setminus F|\geq 2k+1=2|Z'\setminus Z|+1$, and, hence, (ii).
Clearly, (iii) holds in both cases.

Next, assume that $H$ has Type b).
If $k=2$, then let $Z'=Z\cup \{ v_k\}$.
Since $|Z'\setminus Z|=1$, we obtain (i).
Since $p(v_0),p(v_1)\in {\cal F}(Z')\setminus F$ and $Z'\setminus Z\subseteq {\cal F}(Z')\setminus F$, 
we obtain $|{\cal F}(Z')\setminus F|\geq 2|Z'\setminus Z|+1$, and, hence, (ii).
If $k\geq 3$, then let $Z'=Z\cup \{ v_k\}\cup \{ p(v_i):0\leq i\leq k-3\}$.
Since $|Z'\setminus Z|=k-1=n(H)-2\leq 2\log_2(n)$, we obtain (i).
Since $v_1,\ldots,v_{k-1},p(v_{k-2}),p(v_{k-1})\in {\cal F}(Z')\setminus F$ and $Z'\setminus Z\subseteq {\cal F}(Z')\setminus F$, 
we obtain $|{\cal F}(Z')\setminus F|\geq 2k\geq 2|Z'\setminus Z|+1$, and, hence, (ii).
Clearly, (iii) holds in both cases.

Next, assume that $H$ has Type c).
Let $Z'=Z\cup \{ p(v_i):0\leq i\leq k-2\}$.
Since $|Z'\setminus Z|=k-1=n(H)-2\leq 2\log_2(n)$, we obtain (i).
Since $v_1,\ldots, v_{k-1},p(v_{k-1}),p(v_k)\in {\cal F}(Z')\setminus F$ and $Z'\setminus Z\subseteq {\cal F}(Z')\setminus F$, 
we obtain $|{\cal F}(Z')\setminus F|\geq 2k\geq 2|Z'\setminus Z|+1$, and, hence, (ii).
Clearly, (iii) holds.

Next, assume that $H$ has Type d).
Note that, since $G$ has girth at least $5$, we have $\ell\geq 5$.
Let $Z'=Z\cup \{ u_\ell\}\cup \{ p(u_j):2\leq j\leq \ell-2\}$.
Since $|Z'\setminus Z|=\ell-2=n(H)-2\leq 2\log_2(n)$, we obtain (i).
Since $u_2,\ldots, u_{\ell-1},p(u_{\ell-1}),p(u_\ell)\in {\cal F}(Z')\setminus F$ and $Z'\setminus Z\subseteq {\cal F}(Z')\setminus F$, 
we obtain $|{\cal F}(Z')\setminus F|\geq 2\ell-2\geq 2|Z'\setminus Z|+1$, and, hence, (ii).
Clearly, (iii) holds.

Finally, assume that $H$ has Type e).
Let $Z'=Z\cup \{ p(v_i):0\leq i\leq k-1\}\cup \{ u_\ell \}\cup \{ p(u_j):2\leq j\leq \ell-2\}$.
Since $|Z'\setminus Z|=k+\ell-2=n(H)-2\leq 2\log_2(n)$, we obtain (i).
Since $v_1,\ldots,v_k,u_2\ldots,u_{\ell-1},p(u_{\ell-1}),p(u_\ell)\in {\cal F}(Z')\setminus F$ and $Z'\setminus Z\subseteq {\cal F}(Z')\setminus F$, 
we obtain $|{\cal F}(Z')\setminus F|\geq 2k+2\ell-2\geq 2|Z'\setminus Z|+1$, and, hence, (ii).
Clearly, (iii) holds.

This completes the proof of the claim. $\Box$

\medskip

\noindent Since $G$ has maximum degree $3$, we have $n\geq 4$, 
which implies $\frac{1}{2}-\frac{1}{8\log_2(n)+2}\geq \frac{4}{9}$.
For some vertex $v$ of degree $3$, and some neighbor $u$ of $v$,
let $Z_0=N_G[v]\setminus \{ u\}$.
Since $|Z_0|=3$ and $|{\cal F}(Z_0)|\geq 4$, 
we obtain
$$\frac{|Z_0|-2}{|{\cal F}(Z_0)|}\leq \frac{1}{2}-\frac{1}{8\log_2(n)+2}.$$
Clearly, ${\cal F}(Z_0)$ induces a connected subgraph of $G$ of order at least $3$.

Suppose that $Z$ is a set of vertices of $G$ that satisfies the hypotheses of Claim \ref{claim1} such that 
\begin{eqnarray}\label{ez}
\frac{|Z|-2}{|{\cal F}(Z)|}\leq \frac{1}{2}-\frac{1}{8\log_2(n)+2}.
\end{eqnarray}
By Claim \ref{claim1}, the set $Z$ can be extended to a set $Z'$ with the properties stated in Claim \ref{claim1}.
In particular, 
$$\frac{|Z'\setminus Z|}{|{\cal F}(Z')\setminus {\cal F}(Z)|}
\stackrel{(ii)}{\leq} \frac{|Z'\setminus Z|}{2|Z'\setminus Z|+1}
\stackrel{(i)}{\leq} \frac{2\log_2(n)}{4\log_2(n)+1}
=\frac{1}{2}-\frac{1}{8\log_2(n)+2},$$
which implies
$$\frac{|Z'|-2}{|{\cal F}(Z')|}
=\frac{(|Z|-2)+|Z'\setminus Z|}{|{\cal F}(Z)|+|{\cal F}(Z')\setminus {\cal F}(Z)|}
\leq \frac{1}{2}-\frac{1}{8\log_2(n)+2}.$$
In view of the set $Z_0$ defined above, 
this implies the existence of a set $Z$ of vertices of $G$ that satisfies (\ref{ez}) 
such that $F={\cal F}(Z)$ induces a connected subgraph of $G$ of order at least $3$,
and all vertices in $R=V(G)\setminus F$ have degree $1$.
Since $G$ is connected, and every vertex in $F$ has at most two neighbors in $R$, we obtain $|R|\leq 2|F|$.
Since $n=|F|+|R|$, this implies $|F|\geq \frac{n}{3}$ and $|R|\leq \frac{2n}{3}$.
Note that every vertex $v$ in $F$ that has a neighbor in $R$ has exactly two neighbors in $R$.
Let $\tilde{Z}$ arise from $Z$ by adding, for every such vertex $v$ in $F$, 
exactly one of its two neighbors in $R$ to $Z$. 
Clearly, $\tilde{Z}$ is a zero forcing set of $G$, and we obtain
\begin{eqnarray*}
|\tilde{Z}|-2 & = & (|Z|-2)+\frac{1}{2}|R|\\
& \stackrel{(\ref{ez})}{\leq} & \left(\frac{1}{2}-\frac{1}{8\log_2(n)+2}\right)|F|+\frac{1}{2}|R|\\
& \leq & \left(\left(\frac{1}{2}-\frac{1}{8\log_2(n)+2}\right)\cdot \frac{1}{3}+\frac{1}{2}\cdot \frac{2}{3}\right)n\\
& = & \left(\frac{1}{2}-\frac{1}{24\log_2(n)+6}\right)n,
\end{eqnarray*}
which completes the proof. $\Box$

\medskip

\noindent We proceed to our probabilistic upper bound.
For a set $N$ and a non-negative integer $i$, let ${N\choose i}$ be the set of all subsets of $N$ of order $i$.

\begin{theorem}\label{theorem2}
If $G$ is a graph, then
$$Z(G)\leq \sum\limits_{u\in V(G)}\sum\limits_{i=0}^{d_G(u)}(-1)^i\sum\limits_{I\in {N_G(u)\choose i}}\left|\{ u\}\cup\bigcup\limits_{v\in I}N_G[v]\right|^{-1}.$$
\end{theorem}
{\it Proof:} Let $u_1,\ldots,u_n$ be a linear order of the vertices of $G$ selected uniformly at random.
Let $Z$ be the set of those vertices $u_i$
such that $u_i$ is not the unique neighbor within $\{ u_i,\ldots,u_n\}$
of some vertex $u_j$ with $j<i$.
Clearly, $Z$ is a zero forcing set of $G$.
Hence, by the first moment method, $Z(G)\leq \mathbb{E}[|Z|]$.

Let $u$ be a vertex of $G$.
For $v\in N_G(u)$, let $A_v$ be the event that $u$ is the rightmost vertex from $N_G[v]$ within the linear order $u_1,\ldots,u_n$,
that is, if $u=u_j$, then $i<j$ for every $i$ in $\{ 1,\ldots,n\}$ with $u_i\in N_G[v]\setminus \{ u\}$. 
The definition of $Z$ implies
$$\mathbb{P}[u\in Z] = \mathbb{P}\left[\overline{\bigcup\limits_{v\in N_G(u)}A_v}\right].$$
Let $N=\{ u\}\cup\bigcup\limits_{v\in N_G(u)}N_G[v]$ and $d=|N|$.
Note that there are $d!$ linear orders of $N$.
Furthermore, if $I$ is a subset of $N_G(u)$, 
then the number of linear orders $\sigma$ of $N$ 
such that $u$ is the rightmost vertex from $\{ u\}\cup\bigcup\limits_{v\in I}N_G[v]$
within $\sigma$ is exactly $\frac{d!}{\left|\{ u\}\cup\bigcup\limits_{v\in I}N_G[v]\right|}$,
which implies
$$\mathbb{P}\left[\bigcap\limits_{v\in I}A_v\right]=\left|\{ u\}\cup\bigcup\limits_{v\in I}N_G[v]\right|^{-1}.$$
By inclusion-exclusion, we obtain 
\begin{eqnarray*}
\mathbb{P}[u\in Z] & = & \mathbb{P}\left[\overline{\bigcup\limits_{v\in N_G(u)}A_v}\right]\\
&=& \sum\limits_{i=0}^{d_G(u)}(-1)^i\sum\limits_{I\in {N_G(u)\choose i}}\mathbb{P}\left[\bigcap\limits_{v\in I}A_v\right]\\
&=& \sum\limits_{i=0}^{d_G(u)}(-1)^i\sum\limits_{I\in {N_G(u)\choose i}}\left|\{ u\}\cup\bigcup\limits_{v\in I}N_G[v]\right|^{-1}.
\end{eqnarray*}
By linearity of expectation, we have $\mathbb{E}[|Z|]=\sum\limits_{u\in V(G)}\mathbb{P}[u\in Z]$,
and the desired result follows. $\Box$

\medskip

\noindent Since the bound in Theorem \ref{theorem2} is not very explicit, we derive some more explicit corollaries.

For a positive integer $r$, let $H_r=\sum\limits_{i=1}^r\frac{1}{i}$.
It is known that $\lim\limits_{r\to \infty}(H_r-\ln r)\approx 0,577$.

\begin{corollary}\label{corollary1}
If $G$ is a $r$-regular graph of order $n$ and girth at least $5$, then
$$Z(G)\leq \left(\prod_{i=1}^r\left(1-\frac{1}{ri+1}\right)\right)n
=\left(1-\frac{H_r}{r}\right)n+O\left(\left(\frac{H_r}{r}\right)^2\right)n.$$
\end{corollary}
{\it Proof:} By Theorem \ref{theorem2}, we obtain
\begin{eqnarray*}
\frac{Z(G)}{n} & \leq & \frac{1}{n}\sum\limits_{u\in V(G)}\sum\limits_{i=0}^{d_G(u)}(-1)^i\sum\limits_{I\in {N_G(u)\choose i}}\left|\{ u\}\cup\bigcup\limits_{v\in I}N_G[v]\right|^{-1}\\
& = & \sum\limits_{i=0}^{r}(-1)^i{r\choose i}\frac{1}{ri+1}\,\,\,\,\,\,\,\,\,\,\,\,(\mbox{using the regularity and the girth condition})\\
& = & \sum\limits_{i=0}^{r}(-1)^i{r\choose i}\int_{0}^1x^{ri}dx\\
& = & \int_{0}^1\sum\limits_{i=0}^{r}(-1)^i{r\choose i}x^{ri}dx\\
& = & \int_{0}^1(1-x^r)^rdx\,\,\,\,\,\,\,\,\,\,\,\,(\mbox{using the binomial theorem})\\
& = & \frac{1}{r}\int_{0}^1(1-z)^rz^{\frac{1}{r}-1}dz\,\,\,\,\,\,\,\,\,\,\,\,(\mbox{substituting $z=x^r$})\\
& = & \frac{1}{r}B\left(r+1,\frac{1}{r}\right)\,\,\,\,\,\,\,\,\,\,\,\,(\mbox{where $B(\cdot,\cdot)$ is the Beta function})\\
& = & \frac{1}{r}\frac{\Gamma(r+1)\Gamma\left(\frac{1}{r}\right)}{\Gamma\left(1+r+\frac{1}{r}\right)}\,\,\,\,\,\,\,\,\,\,\,\,(\mbox{where $\Gamma(\cdot)$ is the Gamma function})\\
& = & \frac{r!}{(r+\frac{1}{r})(r-1+\frac{1}{r})\ldots (1+\frac{1}{r})}\,\,\,\,\,\,\,\,\,\,\,\,(\mbox{using $\Gamma(x+1)=x\Gamma(x)$})\\
& = & \prod_{i=1}^r\frac{i}{i+\frac{1}{r}}\\
& = & \prod_{i=1}^r\left(1-\frac{1}{ri+1}\right),
\end{eqnarray*}
which implies the first stated bound for $Z(G)$.

Note that 
\begin{eqnarray*}
\prod_{i=1}^r\left(1-\frac{1}{ri+1}\right) 
& = & 1-\sum\limits_{i=1}^r\frac{1}{ri+1}+\sum_{i=2}^r(-1)^i\sum_{I\in {[r]\choose i}}\prod_{j\in I}\frac{1}{rj+1}\\
& = & 1-\left(\sum\limits_{i=1}^r\frac{1}{ri}-\sum\limits_{i=1}^r\frac{1}{ri(ri+1)}\right)+\sum_{i=2}^r(-1)^i\sum_{I\in {[r]\choose i}}\prod_{j\in I}\frac{1}{rj+1}\\
& = & 1-\left(\frac{H_r}{r}-\sum\limits_{i=1}^r\frac{1}{ri(ri+1)}\right)+\sum_{i=2}^r(-1)^i\sum_{I\in {[r]\choose i}}\prod_{j\in I}\frac{1}{rj+1}.
\end{eqnarray*}
Since
\begin{eqnarray*}
\left|\sum_{i=1}^r\frac{1}{ri(ri+1)}\right| & \leq & \frac{1}{r^2}\sum_{i=1}^r\frac{1}{i^2}\leq \frac{\pi^2}{6r^2}=O\left(\left(\frac{H_r}{r}\right)^2\right)
\end{eqnarray*}
and
\begin{eqnarray*}
\left|\sum_{i=2}^r(-1)^i\sum_{I\in {[r]\choose i}}\prod_{j\in I}\frac{1}{rj+1}\right| & \leq &
\sum_{i=2}^r\sum_{I\in {[r]\choose i}}\prod_{j\in I}\frac{1}{rj}\\
& = & \sum_{i=2}^r\frac{1}{r^i}\sum_{I\in {[r]\choose i}}\prod_{j\in I}\frac{1}{j}\\
&\leq & \sum_{i=2}^r\frac{1}{r^i}\frac{1}{i!}\left(\sum_{j_1=1}^r\frac{1}{j_1}\right)\left(\sum_{j_2=1}^r\frac{1}{j_2}\right)\ldots \left(\sum_{j_i=1}^r\frac{1}{j_i}\right)\\
&=& \sum_{i=2}^r\frac{1}{i!}\left(\frac{H_r}{r}\right)^i\\
& \leq & \left(\frac{H_r}{r}\right)^2\sum_{i=2}^r\frac{1}{i!}\\
& \leq & e\left(\frac{H_r}{r}\right)^2\\
&=& O\left(\left(\frac{H_r}{r}\right)^2\right),
\end{eqnarray*}
we obtain the second stated bound for $Z(G)$. $\Box$

\medskip

\noindent Note that 
$$\prod_{i=1}^r\left(1-\frac{1}{ri+1}\right)=
\left\{
\begin{array}{ll}
\frac{81}{140}\approx 0.579 &\mbox{, for $r=3$},\\
\frac{2048}{3315}\approx 0.618 &\mbox{, for $r=4$, and }\\
\frac{15625}{24024}\approx 0.65 &\mbox{, for $r=5$.}
\end{array}\right.$$
In fact, this expression is less than the factor $\frac{r-2}{r-1}$ from (\ref{eb2}) for $r\geq 4$.

If $G$ is a cubic triangle-free graph such that no component of $G$ is $K_{3,3}$, then, 
for every vertex $u$ of $G$, 
the subgraph of $G$ 
that contains all vertices at distance at most $2$ from $u$
as well as all edges incident with neighbors of $u$
is of one of the seven types illustrated in Table \ref{tab1}.
This defines the type of the vertex $u$.

\begin{table}[H]
\begin{center}
\begin{tabular}{|c|c|c|c|c|c|c|}\hline
Type 1 
& 
Type 2 
& 
Type 3 
& 
Type 4 
& 
Type 5 
& 
Type 6 
& 
Type 7 
\\ \hline
 &&&&&&\\ 
\unitlength 0.6mm 
\linethickness{0.4pt}
\ifx\plotpoint\undefined\newsavebox{\plotpoint}\fi 
\begin{picture}(30,30)(0,0)
\put(15,5){\circle*{2}}
\put(15,15){\circle*{2}}
\put(5,15){\circle*{2}}
\put(25,15){\circle*{2}}
\put(12,25){\circle*{2}}
\put(22,25){\circle*{2}}
\put(2,25){\circle*{2}}
\put(18,25){\circle*{2}}
\put(28,25){\circle*{2}}
\put(8,25){\circle*{2}}
\multiput(12,25)(.03370787,-.11235955){89}{\line(0,-1){.11235955}}
\multiput(22,25)(.03370787,-.11235955){89}{\line(0,-1){.11235955}}
\multiput(2,25)(.03370787,-.11235955){89}{\line(0,-1){.11235955}}
\multiput(15,15)(.03370787,.11235955){89}{\line(0,1){.11235955}}
\multiput(25,15)(.03370787,.11235955){89}{\line(0,1){.11235955}}
\multiput(5,15)(.03370787,.11235955){89}{\line(0,1){.11235955}}
\put(5,15){\line(1,-1){10}}
\put(15,5){\line(0,1){10}}
\put(25,15){\line(-1,-1){10}}
\put(15,0){\makebox(0,0)[cc]{$u$}}
\end{picture} 
&
\unitlength 0.6mm 
\linethickness{0.4pt}
\ifx\plotpoint\undefined\newsavebox{\plotpoint}\fi 
\begin{picture}(30,30)(0,0)
\put(15,5){\circle*{2}}
\put(15,15){\circle*{2}}
\put(5,15){\circle*{2}}
\put(25,15){\circle*{2}}
\put(22,25){\circle*{2}}
\put(2,25){\circle*{2}}
\put(18,25){\circle*{2}}
\put(28,25){\circle*{2}}
\multiput(22,25)(.03370787,-.11235955){89}{\line(0,-1){.11235955}}
\multiput(2,25)(.03370787,-.11235955){89}{\line(0,-1){.11235955}}
\multiput(15,15)(.03370787,.11235955){89}{\line(0,1){.11235955}}
\multiput(25,15)(.03370787,.11235955){89}{\line(0,1){.11235955}}
\put(5,15){\line(1,-1){10}}
\put(15,5){\line(0,1){10}}
\put(25,15){\line(-1,-1){10}}
\put(15,0){\makebox(0,0)[cc]{$u$}}
\put(10,25){\circle*{2}}
\put(15,15){\line(-1,2){5}}
\put(10,25){\line(-1,-2){5}}
\end{picture}
&
\unitlength 0.6mm 
\linethickness{0.4pt}
\ifx\plotpoint\undefined\newsavebox{\plotpoint}\fi 
\begin{picture}(30,30)(0,0)
\put(15,5){\circle*{2}}
\put(15,15){\circle*{2}}
\put(5,15){\circle*{2}}
\put(25,15){\circle*{2}}
\put(2,25){\circle*{2}}
\put(5,15){\line(1,-1){10}}
\put(15,5){\line(0,1){10}}
\put(25,15){\line(-1,-1){10}}
\put(15,0){\makebox(0,0)[cc]{$u$}}
\put(10,25){\circle*{2}}
\put(10,25){\line(-1,-2){5}}
\put(20,25){\circle*{2}}
\put(28,25){\circle*{2}}
\put(10,25){\line(1,-2){5}}
\put(15,15){\line(1,2){5}}
\put(10,25){\line(3,-2){15}}
\multiput(2,25)(.03370787,-.11235955){89}{\line(0,-1){.11235955}}
\multiput(28,25)(-.03370787,-.11235955){89}{\line(0,-1){.11235955}}
\end{picture}
&
\unitlength 0.6mm 
\linethickness{0.4pt}
\ifx\plotpoint\undefined\newsavebox{\plotpoint}\fi 
\begin{picture}(30,30)(0,0)
\put(15,5){\circle*{2}}
\put(15,15){\circle*{2}}
\put(5,15){\circle*{2}}
\put(25,15){\circle*{2}}
\put(2,25){\circle*{2}}
\put(5,15){\line(1,-1){10}}
\put(15,5){\line(0,1){10}}
\put(25,15){\line(-1,-1){10}}
\put(15,0){\makebox(0,0)[cc]{$u$}}
\put(10,25){\circle*{2}}
\put(10,25){\line(-1,-2){5}}
\put(20,25){\circle*{2}}
\put(28,25){\circle*{2}}
\put(10,25){\line(1,-2){5}}
\multiput(2,25)(.03370787,-.11235955){89}{\line(0,-1){.11235955}}
\multiput(28,25)(-.03370787,-.11235955){89}{\line(0,-1){.11235955}}
\multiput(2,25)(.0437710438,-.0336700337){297}{\line(1,0){.0437710438}}
\put(20,25){\line(1,-2){5}}
\end{picture}
&
\unitlength 0.6mm 
\linethickness{0.4pt}
\ifx\plotpoint\undefined\newsavebox{\plotpoint}\fi 
\begin{picture}(30,30)(0,0)
\put(15,5){\circle*{2}}
\put(15,15){\circle*{2}}
\put(5,15){\circle*{2}}
\put(25,15){\circle*{2}}
\put(2,25){\circle*{2}}
\put(5,15){\line(1,-1){10}}
\put(15,5){\line(0,1){10}}
\put(25,15){\line(-1,-1){10}}
\put(15,0){\makebox(0,0)[cc]{$u$}}
\put(10,25){\circle*{2}}
\put(10,25){\line(-1,-2){5}}
\put(20,25){\circle*{2}}
\put(28,25){\circle*{2}}
\put(10,25){\line(1,-2){5}}
\multiput(2,25)(.03370787,-.11235955){89}{\line(0,-1){.11235955}}
\multiput(28,25)(-.03370787,-.11235955){89}{\line(0,-1){.11235955}}
\put(20,25){\line(1,-2){5}}
\put(20,25){\line(-1,-2){5}}
\end{picture}
&
\unitlength 0.6mm 
\linethickness{0.4pt}
\ifx\plotpoint\undefined\newsavebox{\plotpoint}\fi 
\begin{picture}(21,30)(0,0)
\put(10,5){\circle*{2}}
\put(10,15){\circle*{2}}
\put(10,25){\circle*{2}}
\put(0,15){\circle*{2}}
\put(0,25){\circle*{2}}
\put(20,15){\circle*{2}}
\put(20,25){\circle*{2}}
\put(0,15){\line(1,-1){10}}
\put(10,5){\line(0,1){10}}
\put(20,15){\line(-1,-1){10}}
\put(10,0){\makebox(0,0)[cc]{$u$}}
\put(0,15){\line(0,1){10}}
\put(0,25){\line(1,-1){10}}
\put(10,15){\line(0,1){10}}
\put(10,25){\line(-1,-1){10}}
\put(20,15){\line(0,1){10}}
\put(20,15){\line(-1,1){10}}
\end{picture}
&
\unitlength 0.6mm 
\linethickness{0.4pt}
\ifx\plotpoint\undefined\newsavebox{\plotpoint}\fi 
\begin{picture}(21,30)(0,0)
\put(10,5){\circle*{2}}
\put(10,15){\circle*{2}}
\put(10,25){\circle*{2}}
\put(0,15){\circle*{2}}
\put(0,25){\circle*{2}}
\put(20,15){\circle*{2}}
\put(20,25){\circle*{2}}
\put(0,15){\line(1,-1){10}}
\put(10,5){\line(0,1){10}}
\put(20,15){\line(-1,-1){10}}
\put(10,0){\makebox(0,0)[cc]{$u$}}
\put(10,15){\line(-1,1){10}}
\put(0,25){\line(0,-1){10}}
\put(0,15){\line(1,1){10}}
\put(10,25){\line(1,-1){10}}
\put(20,15){\line(0,1){10}}
\put(20,25){\line(-1,-1){10}}
\end{picture}
\\ 
 &&&&&&\\
\hline
 &&&&&&\\
$p_1=\frac{81}{140}$ 
& 
$p_2=\frac{149}{252}$ 
& 
$p_3=\frac{5}{8}$ 
& 
$p_4=\frac{171}{280}$ 
& 
$p_5=\frac{101}{168}$ 
& 
$p_6=\frac{269}{420}$
& 
$p_7=\frac{17}{28}$\\
 &&&&&&\\ \hline
\end{tabular}
\caption{The seven possible types of the vertex $u$.}\label{tab1}
\end{center}
\end{table}

\begin{corollary}\label{corollary2}
If $G$ is a cubic triangle-free graph such that no component of $G$ is $K_{3,3}$,
and $G$ has $n_i$ vertices of type $i$ for $i\in \{ 1,\ldots,7\}$, then 
$Z(G)\leq \sum\limits_{i=1}^7p_in_i$.
\end{corollary}
{\it Proof:} This follows immediately from Theorem \ref{theorem2} by calculating the probabilities $\mathbb{P}[u\in Z]$
considered within the proof of Theorem \ref{theorem2} for the vertices $u$ of the different types.
If $u$ has type $4$ for instance, then 
$\mathbb{P}[u\in Z]=1-\frac{3}{4}+\frac{1}{5}+\frac{2}{7}-\frac{1}{8}=\frac{171}{280}$. $\Box$

\medskip

\noindent We proceed to the proof of two further cases of the conjecture of Davila and Kenter.

\begin{theorem}\label{theorem1}
If $G$ is a graph of girth $g$ in $\{ 5,6\}$ and minimum degree $\delta$ at least $2$, then $$Z(G)\geq (g-2)(\delta-2)+2.$$
\end{theorem}
{\it Proof:} Let $G$ be as in the statement.
Let $Z$ be a zero forcing set of minimum cardinality. For a contradiction, suppose that $|Z|\leq (g-2)(\delta-2)+1$. 
For $\delta=2$, this implies that $Z$ contains exactly one vertex, say $v_1$.
Since $G$ has more than one vertex,
and $v_1$ has degree at least $2$,
no vertex in $V(G)\setminus Z$ is the unique neighbor of $v_1$,
which implies a contradiction.
Hence, $\delta\geq 3$.
Since $g\geq 5$, the order $n$ of $G$ is at least $1+\delta+\delta(\delta-1)=\delta^2+1$.
Since $g\in \{ 5,6\}$ and $\delta\geq 3$,
we obtain $n-|Z|\geq \delta^2-(g-2)(\delta-2)\geq g-2$,
which implies that a forcing sequence
${\cal S}:v_1\to u_1, v_2\to u_2,\ldots,v_k\to u_k$
satisfies $k\geq g-2$.
Let $Z'=\{ v_1,\ldots,v_{g-2}\}$.
Let $N=\left(\bigcup\limits_{v\in Z'}N_G(v)\right)\setminus Z'$.
Since ${\cal S}$ is a forcing sequence, 
$Z'\cup N\subseteq Z\cup \{ u_1,\ldots,u_{g-2}\}$, 
and, hence, 
\begin{eqnarray*}
|N| & = & |Z'\cup N|-(g-2)\\
& \leq & |Z\cup \{ u_1,\ldots,u_{g-2}\}|-(g-2)\\
& = & |Z|\\
& \leq & (g-2)(\delta-2)+1.
\end{eqnarray*}
Let $G'=G[Z']$.
Let $G'$ have $\kappa$ components.
Note that $\kappa\leq |Z'|=g-2\leq 4$.
Let $m'$ be the number of edges between $Z'$ and $N$.
Since $g>|Z'|$, the graph $G'$ is a forest and has exactly $g-\kappa-2$ edges.
This implies that 
$m'\geq (g-2)\delta-2(g-\kappa-2)\geq |N|+2\kappa-1$.
Since $m'>|N|$, some vertex in $N$ has more than one neighbor in $Z'$.
Since $g>|Z'|+1$, 
no vertex in $N$ has two neighbors in the same component of $G'$.
This implies that $\kappa\geq 2$.

First, we assume that $\kappa=2$.
If three vertices in $N$ have neighbors in both components of $G'$,
then $G$ has a cycle of length at most $g-1$, which is a contradiction.
Hence, at most two neighbors in $N$ have neighbors in both components of $G'$, which implies the contradiction
$m'\leq |N|+2<|N|+2\kappa-1$.

Next, we assume that $\kappa=3$.
If some vertex $u$ in $N$ has neighbors in all three components of $G'$ and another vertex $u'$ has two neighbors in $Z'$,
then $G$ has a cycle of length at most $g-1$, which is a contradiction.
Similarly, 
if two distinct vertices in $N$ have neighbors in the same two components of $G'$, then $G$ has a cycle of length at most $g-1$, which is a contradiction.
These observations imply the contradiction
$m'\leq |N|+{3\choose 2}<|N|+2\kappa-1$.

Finally, we assume that $\kappa=4$, which implies that $g=6$, and that $Z'$ is an independent set.
Again, no two distinct vertices in $N$ have neighbors in the same two components of $G'$.
If some vertex in $N$ has neighbors in three components of $G'$, then this implies the contradiction
$m'\leq |N|+5<|N|+2\kappa-1$.
Similarly, if no vertex in $N$ has neighbors in three components of $G'$, then this implies the contradiction
$m'\leq |N|+{4\choose 2}<|N|+2\kappa-1$.

This final contradiction completes the proof. $\Box$

\medskip

\noindent {\bf Acknowledgment} We thank Henning Bruhn-Fujimoto for fruitful discussion.

\end{document}